\documentclass[10pt,reqno]{amsart}

\usepackage{amsmath,amsfonts,amsthm,amssymb}
\usepackage{latexsym}
\usepackage{enumerate}
\usepackage[mathscr]{eucal}
\usepackage{eqlist}
\usepackage{geometry}

\makeatletter
\@namedef{subjclassname@2010}{
  \textup{2010} Mathematics Subject Classification}
\makeatother

\newtheorem{theorem}{Theorem}[section]
\newtheorem{proposition}{Proposition}[section]
\newtheorem{lemma}{Lemma}[section]
\newtheorem{corollary}{Corollary}[section]

\theoremstyle{definition}

\theoremstyle{remark}
\newtheorem{remark}{Remark}[section]
\newtheorem*{remarks}{Remarks}

\numberwithin{equation}{section}

\title
{Three-Dimensional Almost Kenmotsu Manifolds \\ Satisfying Certain Nullity Conditions}

\author{Vincenzo Saltarelli}  

\address{Department of Mathematics, University of Bari \newline\indent
         Via E. Orabona, 4 \newline\indent
         I--70125 Bari, ITALY}

\email{saltarelli@dm.uniba.it}



\date{}

\begin{document}

\begin{abstract}
We investigate $3$-dimensional almost Kenmotsu manifolds satisfying special types of  nullity conditions depending on two smooth functions $\kappa,\mu$. When either $\kappa<-1$ and $\mu=0$ or $h=0$, such conditions coincide with the $\kappa$-nullity condition which we show to be equivalent to the $\eta$-Einstein one. As an application of this result, we obtain examples of $N(\kappa)$-quasi Einstein manifolds.
Moreover, for the aforementioned manifolds, some complete local descriptions of their structure are established, building local ``models'' for each of them.
\end{abstract}

\maketitle

\section{Introduction}

The so-called nullity conditions on different kinds of almost contact metric manifolds have known an increasing interest in the last 15 years and they were widely studied in various recent papers.
In \cite{BKP} D. E. Blair,  T. Koufogiorgos and B. J. Papantoniou introduced  a nullity condition, called $(\kappa,\mu)$-nullity condition with $\kappa,\mu\in\mathbb{R}$, on manifolds endowed with a contact metric structure $(\varphi,\xi,\eta,g)$. Afterwards, a natural generalization of this condition has been considered in \cite{GX1}, \cite{GX2}, \cite{KT}, allowing $\kappa,\mu$ to be smooth functions. This means that the curvature of the Levi-Civita connection satisfies
\begin{equation}\label{nullitycond1}
R(X,Y)\xi=\kappa\big(\eta(Y)X-\eta(X)Y\big)+\mu\big(\eta(Y)hX-\eta(X)hY\big),
\end{equation}
for all vector fields $X,Y$ and for some $\kappa,\mu$ in the ring $\mathfrak{F}(M)$ of the smooth functions on $M$, where $h:=(1/2)\mathcal{L}_{\xi}\varphi$, $\mathcal{L}$ denoting the Lie differentiation. In particular, in \cite{KT} it is proved that only in dimension $3$ there are examples of contact metric manifold satisfying (\ref{nullitycond1}) with $\kappa,\mu$ non-constant functions. This problematic has been studied by P. Dacko and Z. Olszak \cite{DO2} in the framework of almost cosymplectic manifolds as well. In this paper we deal with almost Kenmotsu manifolds, focusing mainly on the $3$-dimensional non-normal case. Normal almost Kenmotsu manifolds, known as Kenmotsu manifolds, are closely related with warped products of two Riemannian manifolds, which set up one of the classes of the Tanno's classification \cite{Ta1}. In general, the structure of normal almost contact metric manifolds of dimension $3$ have been recently further studied in \cite{DK}.

Almost Kenmotsu manifolds satisfying (\ref{nullitycond1}) with $\kappa,\mu\in\mathbb{R}$ were studied in \cite{DP2}, proving that $h=0$ and $\kappa=-1$. For this reason, a modified nullity condition involving the tensor $h'=h\circ\varphi$ was introduced and studied in the same paper. It has been later generalized in \cite{PS2}, requiring that $\xi$ belongs to the generalized $(\kappa,\mu)'$-nullity distribution, i.e. for all vector fields $X,Y$
\begin{equation}\label{nullitycond2}
R(X,Y)\xi=\kappa\big(\eta(Y)X-\eta(X)Y\big)+\mu\big(\eta(Y)h'X-\eta(X)h'Y\big),
\end{equation}
for some $\kappa,\mu\in\mathfrak{F}(M)$. If an almost Kenmotsu manifold satisfies (\ref{nullitycond2}), then $\kappa\leq-1$ \cite{PS2}. When $\kappa=-1$ identically, then $h'=0$ and $M$ is locally a warped product of an almost K\"{a}hler manifold and an open interval. In the special case $\kappa,\mu\in\mathbb{R}$, it was proved in \cite{DP2} that, if $\kappa<-1$, then $\mu=-2$ and $h'$ admits $3$ eigenvalues $\lambda, -\lambda, 0$, with $0$ as simple eigenvalue and $\lambda=\sqrt{-1-\kappa}$.
Furthermore, the following classification theorem was established:

\begin{theorem}[\cite{DP2}]\label{classification theorem}
Let $(M^{2n+1},\varphi,\xi,\eta,g)$ be an almost Kenmotsu manifold such that $\xi$ belongs to the $(\kappa,-2)'$-nullity distribution and $h'\neq0$. Then, $M^{2n+1}$ is locally isometric to the warped product $\mathbb{H}^{n+1}(\kappa-2\lambda)\times_{f}\mathbb{R}^n$ or $B^{n+1}(\kappa+2\lambda)\times_{f'}\mathbb{R}^n$, where $f=ce^{(1-\lambda)t}$ and $f'=c'e^{(1+\lambda)t}$, with $c,c'$ positive constants.
\end{theorem}

In \cite{PS2} A. M. Pastore and the author proved that, in dimension greater than or equal to $5$, the functions $\kappa$ and $\mu$ in (\ref{nullitycond1}) or (\ref{nullitycond2}) can vary only in the direction of $\xi$
(i.e. $d\kappa\wedge\eta=0$) and the curvature is completely determined if $\kappa<-1$. In general, similar results are not true in dimension $3$.
However, some explicit examples of almost Kenmotsu manifolds satisfying (\ref{nullitycond1}) and (\ref{nullitycond2}) with non-constant $\kappa,\mu$ were given in any dimension, thus showing a different behavior with respect to the contact case.

The aim of the present paper is to carry on a deeper study of the three dimensional case.
After a section containing basic information about almost Kenmotsu manifolds,
we will start by establishing some useful properties in section \ref{sec3}. In particular, motivated by the growing importance of $N(\kappa)$-quasi Einstein manifolds in the last years (e.g. \cite{TH, YDA}), we prove that $3$-dimensional almost Kenmotsu manifolds whose characteristic vector field $\xi$ belongs to a $\kappa$-nullity distribution are $N(\kappa)$-quasi Einstein manifolds. With section \ref{sec4} we analyze a special class of almost Kenmotsu $3$-manifolds satisfying (\ref{nullitycond1}), which are not Kenmotsu, and we present one of the main results of the paper. We show that the structure of a $3$-dimensional almost Kenmotsu generalized $(\kappa,\mu)$-manifold with $h\neq0$ is completely locally described provided that $d\kappa\wedge\eta=0$ (Theorem 4.1). Section \ref{sec5} deals an analogous situation for the class of almost Kenmotsu generalized $(\kappa,\mu)'$-manifolds of dimension $3$  with $\kappa<-1$. A complete local description of their structure is obtained (Theorem 5.1).
In the last section, with a different approach, we give other important characterizations (Theorem 6.1 and Theorem 6.2) and explicit local models for the two classes considered.

All manifolds are assumed to be smooth (i.e. differentiable of class $\mathcal{C}^{\infty}$) and connected.

\section{Preliminaries}

In this section, we recall some basic data about almost Kenmotsu manifolds and the main properties obtained in \cite{PS2} for almost Kenmotsu manifolds satisfying the generalized nullity conditions.

\subsection{Almost Kenmotsu manifolds}

An almost contact metric manifold is a ($2n+1$)-manifold $M^{2n+1}$ endowed with a structure $(\varphi,\xi,\eta,g)$ given by a $(1,1)$-tensor field $\varphi$, a vector field $\xi$, a $1$-form $\eta$ and a Riemannian metric $g$ satisfying
\begin{align*}
&\varphi^2=-I+\eta\otimes\xi, \quad \eta(\xi)=1,  \\
&g(\varphi X,\varphi Y)=g(X,Y)-\eta(X)\eta(Y) \qquad \forall X,Y\in\Gamma(TM).
\end{align*}
The fundamental $2$-form $\Phi$ associate with the structure is defined by $\Phi(X,Y)=g(X,\varphi Y)$ for any vector fields $X$ and $Y$. The structure is normal if the tensor field  $N=[\varphi,\varphi]+2d\eta\otimes\xi$ vanishes, where $[\varphi,\varphi]$ is the Nijenhuis torsion of $\varphi$. For more details, we refer the reader to \cite{Bl1}.

According to \cite{JV}, an almost contact metric manifold $(M^{2n+1},\varphi,\xi,\eta,g)$ is said to be \emph{almost Kenmotsu manifold} if
\begin{equation}\label{almKen}
d\eta = 0, \qquad d\Phi = 2\eta \wedge \Phi.
\end{equation}
A normal almost Kenmotsu manifold is a Kenmotsu manifold (cf. \cite{JV, Ken}). \\
Let $(M^{2n+1},\varphi,\xi,\eta,g)$ be an almost Kenmotsu manifold. Since $d\eta=0$, the canonical distribution $\mathcal{D}=\ker(\eta)$ orthogonal to $\xi$ is integrable. Denote by $\mathcal{F}$ the $1$-codimensional canonical foliation of $M$ generated by $\mathcal{D}$.
The tensor field $\varphi$ induces an almost complex structure $J$ on any leaf $M'$ of $\mathcal{F}$ and, if $G$ is the Riemannian metric induced by $g$ on $M'$, the pair $(J,G)$ defines an almost Hermitian structure on $M'$. From (\ref{almKen}) we infer that $(J,G)$ is an almost K\"ahlerian structure on $M'$. Furthermore, we have ${\mathcal L}_{\xi} \eta =0$ and $[\xi,X] \in \mathcal D$ for any $X \in \mathcal D$. The Levi-Civita connection fulfills the following relation (cf. \cite{Bl1}):
\begin{equation}\label{Levi-Civita}
2g\big((\nabla_X \varphi)Y,Z\big)=2\eta(Z)g(\varphi X, Y) -2\eta(Y)g(\varphi X,Z)
+ g(N(Y,Z),\varphi X),
\end{equation}
for any vector fields $X,Y,Z$, from which we deduce that $\nabla_{\xi} \varphi =0$, so that $\nabla_{\xi} \xi =0$ and $\nabla_\xi X\in\mathcal D$ for any $X\in \mathcal D$.
From Lemma 2.2 in \cite{KP} we also have, for any $X,Y\in\Gamma(TM)$,
\begin{equation}\label{Levi-Civita bis}
(\nabla_X \varphi)Y+(\nabla_{\varphi X} \varphi)\varphi Y=-2g(X,\varphi  Y)\xi-\eta(Y)\varphi X-\eta(Y)hX,
\end{equation}
where $h:=(1/2)\mathcal{L}_{\xi}\varphi$.

The $(1,1)$-tensor fields $h$ and $h':=h\circ\varphi$ are both self-adjoint operators such that $h(\xi)=h'(\xi)=0$ and satisfy
$$
\eta\circ h=\eta\circ h'=0, \quad h\circ\varphi+\varphi\circ h=0, \quad h'\circ\varphi+\varphi\circ h'=0, \quad \mathrm{tr}(h)=\mathrm{tr}(h')=0.
$$
It follows that non-vanishing $h$ and $h'$ have the non-zero eigenvalues with opposite sign; moreover, they admit the same eigenvalues, but different eigenspaces. If $\lambda\neq0$, we will denote by $[\lambda]$ and $[\lambda]'$ the distributions of eigenvectors of $h$ and $h'$, respectively, with eigenvalue $\lambda$. Moreover, for any vector field $X$, (\ref{Levi-Civita}) implies the following relation
\begin{equation}\label{e:nablaxi}
\nabla_X \xi=X-\eta(X)\xi-\varphi hX.
\end{equation}
In \cite{KP} it is proved that the integral submanifolds of $\mathcal{D}$ are totally umbilical submanifolds of $M^{2n+1}$ if and only if $h=0$, which is equivalent to the vanishing of $h'$. In this case the manifold is locally a warped product $M'\times_{f} N^{2n}$, where $N^{2n}$ is an almost K\"{a}hler manifold, $M'$ is an open interval with coordinate $t$, and $f(t)=ce^{t}$ for some positive constant $c$ (see \cite{DP1}). If, in addition, the integral submanifolds of $\mathcal{D}$ are K\"{a}hler, then $M^{2n+1}$ is a Kenmotsu manifold. In particular, an almost Kenmotsu $3$-manifold with $h=0$ is a Kenmotsu manifold.

Moreover, we recall that the Riemannian curvature $R$ of an almost Kenmotsu manifold satisfies the following general relations (cf. \cite{DP1, DP2}):
\begin{equation}\label{curv1}
R(Y, Z)\xi = \eta(Y)(Z-\varphi hZ)-\eta(Z)(Y-\varphi hY)+(\nabla_Z\varphi h)Y-(\nabla_Y\varphi h)Z
\end{equation}

\begin{equation}\label{l}
\varphi l \varphi-l=2(-\varphi^2  +h^2)
\end{equation}
where $l$ is the self-adjoint operator defined by $l(X):=R(X, \xi)\xi$, for any vector field $X$. The above relations can be also written in terms of $h'$ since $\varphi\circ h=-h'$, $h=\varphi\circ h'$ and $h^2=h'^2$.

Finally, we recall that an almost contact metric manifold $(M^{2n+1},\varphi,\xi,\eta,g)$
is said to be \emph{$\eta$-Einstein} if its Ricci tensor satisfies
\begin{equation*}\label{ricci}
\mathrm{Ric} = a g+ b\eta\otimes\eta
\end{equation*}
or equivalently
\begin{equation}\label{ricci1}
Q = a I+ b\eta\otimes\xi\,,
\end{equation}
where $a, b$ are smooth functions on $M^{2n+1}$. An its extension to every Riemannian manifold is the notion of quasi Einstein manifold.
A \emph{quasi Einstein manifold} in the sense of \cite{CM} is a non-flat Riemannian manifold $(M^n,g)$ whose Ricci operator $Q$ is not identically zero and satisfies (\ref{ricci1}) with $b\neq0$ for a suitable nowhere vanishing $1$-form $\eta$ and unit vector field $\xi$ such that $\eta=g(\cdot,\xi)$. Thus, in particular, any $\eta$-Einstein almost Kenmotsu manifold is quasi Einstein.

\subsection{Nullity distributions}

Let $(M^{2n+1},\varphi,\xi,\eta,g)$ be an almost contact metric manifold and $\kappa,\mu\in\mathfrak{F}(M)$. The \emph{generalized $(\kappa,\mu)$-nullity distribution $N(\kappa,\mu)$} is the distribution defined by putting for each $p\in M^{2n+1}$
\begin{align*}
N_p(\kappa,\mu)= \{Z \in T_p M^{2n+1} \mid R(X, Y)Z&=
\kappa(g(Y,Z)X-g(X,Z)Y)\\
&+ \mu(g(Y,Z)h X-g(X,Z)h Y)\}\,,
\end{align*}
where $X$ and $Y$ are arbitrary vectors in  $T_p M^{2n+1}$. \\
The \emph{generalized $(\kappa,\mu)'$-nullity distribution $N(\kappa,\mu)'$} is obtained by replacing $h$ with $h'$. If $\mu=0$ or $h=0$, both distributions reduce to the well-known $\kappa$-nullity distribution $N(\kappa)$. The generalized $(\kappa,\mu)$-nullity condition (\ref{nullitycond1}) (resp. the generalized $(\kappa,\mu)'$-nullity condition (\ref{nullitycond2})) is obtained by requiring that $\xi$ belongs to some $N(\kappa,\mu)$ (resp. $N(\kappa,\mu)'$).
For convenience, an almost contact metric manifold satisfying (\ref{nullitycond1}) (resp. (\ref{nullitycond2})) will be called \emph{generalized $(\kappa,\mu)$-manifold} (resp. \emph{generalized $(\kappa,\mu)'$-manifold}).

We observe that, in an almost Kenmotsu manifold, if $\xi \in N(\kappa,\mu)$ or $\xi \in N(\kappa,\mu)'$, (\ref{curv1}) implies
\begin{equation} \label{eq:h' codazzi}
(\nabla_X h')Y - (\nabla_Y h')X = 0,
\end{equation}
for any $X,Y \in \mathcal{D}$. Furthermore, in \cite{PS2} the following relations are found:
\begin{equation} \label{eq:h^2}
h^2=h'^2=(\kappa+1)\varphi^2, \quad\quad Q(\xi)=2nk\xi,
\end{equation}
$Q$ being the Ricci operator.
It follows that at every point of an almost Kenmotsu manifold:
\begin{enumerate}[a)]
\item $\kappa\leq-1$;
\item $\kappa=-1$ if and only if $h=0$ or, equivalently, $h'=0$;
\item if $\kappa<-1$, then the eigenvalues of $h$ and $h'$ are $0$ of multiplicity $1$ and $\lambda=\sqrt{-1-\kappa}$ and $-\lambda$, both of them with multiplicity $n$.
\end{enumerate}
In the case of the $(\kappa,\mu)$-nullity condition we also have
\begin{equation}\label{eq:nh}
\nabla_\xi h=-2h -\mu\varphi h,  \qquad d\lambda(\xi)=-2\lambda, \qquad d\kappa(\xi)=-4(\kappa+1).
\end{equation}
Whereas, the belonging of $\xi$ to the $(\kappa,\mu)'$-nullity distribution yields
\begin{equation}\label{eq:nh'}
\nabla_\xi h'=-(\mu+2)h',  \qquad d\lambda(\xi)=-\lambda(\mu+2), \qquad d\kappa(\xi)=-2(\kappa+1)(\mu+2).
\end{equation}

Here we prove the following additional result.

\begin{proposition}
Let $(M^{2n+1},\varphi,\xi,\eta,g)$ be an almost Kenmotsu manifold. If $\xi$ belongs to the generalized $(\kappa,\mu)$-nullity distribution, then one has:
\begin{equation}\label{Lie derivatives1}
\mathcal{L}_{\xi}h=2\lambda^2\varphi-2h+\mu h',\quad \mathcal{L}_{\xi}h'=-\mu h-2h'.
\end{equation}
If $\xi$ belongs to the generalized $(\kappa,\mu)'$-nullity distribution, then one has:
\begin{equation}\label{Lie derivatives2}
\mathcal{L}_{\xi}h'=-(\mu+2)h',\quad \mathcal{L}_{\xi}h=2\lambda^2\varphi-(\mu+2)h.
\end{equation}
\end{proposition}
\begin{proof}
We remind that, given a $(1,1)$-tensor field $T$ on $M$, the following general relation holds
$$
\mathcal{L}_{X}T=\nabla_X T+T\circ\nabla X-(\nabla X)\circ T,
$$
for any $X\in\Gamma(TM)$. Applying it to $T=h$ and $X=\xi$, the required relation follows from the first equation in (\ref{eq:nh}). Now, using the equation just proved and the first relation of (\ref{eq:h^2}) in
$$
\mathcal{L}_{\xi}h'=\mathcal{L}_{\xi}(h\circ\varphi)=(\mathcal{L}_{\xi}h)\varphi + h(\mathcal{L}_{\xi}\varphi)
$$
we get the second formula in (\ref{Lie derivatives1}). Analogously, taking into account that $h=\varphi\circ h'$, the relations in (\ref{Lie derivatives2}) are obtained.
\end{proof}

\section{Some further properties}\label{sec3}

We first establish general formulas which hold on every almost Kenmotsu manifold, without any restriction on the dimension.

\begin{lemma}
Let $(M^{2n+1},\varphi,\xi,\eta,g)$ be an almost Kenmotsu manifold. Then, for any orthonormal frame $\{X_i\}_{1\leq i\leq 2n+1}$, one has
\begin{align}
\sum_{i=1}^{2n+1}(\nabla_{X_i}h')X_i&=Q\xi+2n\xi \label{tracenh'} \\
\sum_{i=1}^{2n+1}(\nabla_{X_i}\varphi)X_i&=0 \label{tracenfi}.
\end{align}
\end{lemma}
\begin{proof}
Let $\{X_i\}_{1\leq i\leq 2n+1}$ be an orthonormal frame. For any vector field $X$, putting $Y=X_i$, replacing $Z$ by $\varphi X$ in (\ref{curv1}) and taking the inner product with $X_i$, we get
\begin{equation}\label{tracenh'0}
\begin{split}
\sum_{i=1}^{2n+1} g(R(X_i, \varphi X)\xi, X_i)&=\sum_{i=1}^{2n+1}\eta(X_i)g(\varphi X-hX, X_i)
                                                  +\sum_{i=1}^{2n+1}g((\nabla_{X_i} h')\varphi X, X_i) \\
                                               &\quad+\sum_{i=1}^{2n+1}g((\nabla_{\varphi X} h')X_i, X_i).
\end{split}
\end{equation}
By definition of the Ricci tensor, the left-hand side is equal to $g(Q\xi, \varphi X)$; while the first term on the right-hand side vanishes, since $\sum_{i=1}^{2n+1}\eta(X_i)X_i=\xi$ and $h\xi=\varphi\xi=0$. Since $\mathrm{tr}(h')=0$, the last term vanishes as well.
Therefore, using the symmetry of $\nabla_{X_i} h'$ and the skew-symmetry of $\varphi$, (\ref{tracenh'0}) reduces to
$$
\sum_{i=1}^{2n+1}\varphi(\nabla_{X_i}h')X_i=\varphi Q\xi.
$$
Applying $\varphi$ to this equation, using $\varphi^2=-I+\eta\otimes\xi$ and taking into account that, by definition, $g(Q\xi,\xi)=Ric(\xi,\xi)=\mathrm{tr}\ l$, we get
\begin{equation} \label{tracenh'1}
\sum_{i=1}^{2n+1}(\nabla_{X_i}h')X_i-\sum_{i=1}^{2n+1}\eta((\nabla_{X_i}h')X_i)\xi=Q\xi-(\mathrm{tr}\ l)\xi.
\end{equation}
Now, using (\ref{e:nablaxi}) and $\mathrm{tr}(h')=0$, one has
\begin{align*}
\sum_{i=1}^{2n+1}\eta((\nabla_{X_i}h')X_i)\xi & = -\sum_{i=1}^{2n+1}g(X_i, h'(\nabla_{X_i}\xi))\xi
= -\sum_{i=1}^{2n+1}g(X_i,h'X_i+h'^2 X_i)\xi  \\
&=-\mathrm{tr}(h')\xi-\mathrm{tr}(h'^2)\xi = -\mathrm{tr}(h'^2)\xi.
\end{align*}
On the other hand, (\ref{l}) implies $K(\varphi X,\xi)+K(X,\xi)=-2-2g(h^2 X,X)$,
from which it follows that $Ric(\xi,\xi)=-2n-\mathrm{tr}(h^2)$.
Hence $\mathrm{tr}(h'^2)=\mathrm{tr}(h^2)=-2n-\mathrm{tr}\ l$ and so, from (\ref{tracenh'1}), we get (\ref{tracenh'}). In order to obtain (\ref{tracenfi}), we point out that the left-hand side is independent of the particular choice of the orthonormal frame. We therefore may compute it by choosing a $\varphi$-basis
$\{E_{i},\varphi E_{i}, \xi\}_{1\leq i\leq n}$. By using (\ref{Levi-Civita bis}) with $X=Y=E_{i}$, since $\eta(E_{i})=0$, we have
\begin{align*}
\sum_{i=1}^{n}(\nabla_{E_i}\varphi)E_i+\sum_{i=1}^{n}(\nabla_{\varphi E_i}\varphi)\varphi E_i=-\sum_{i=1}^{n}\eta(E_i)(\varphi E_i+h E_i)=0,
\end{align*}
from which, being $\nabla_{\xi}\varphi=0$, (\ref{tracenfi}) follows.
\end{proof}

The next lemma concerns almost Kenmotsu manifolds having the canonical distribution $\mathcal{D}$ with K\"{a}hler leaves for which the following formula holds (cf. \cite{FP}):
\begin{equation}\label{kaehler leaves bis}
(\nabla_X\varphi)Y=g(\varphi X+hX,Y)\xi-\eta(Y)(\varphi X+hX)\,, \:\forall X,Y \in \Gamma(TM).
\end{equation}

\begin{lemma}\label{tracenh}
Let $(M^{2n+1},\varphi,\xi,\eta,g)$ be an almost Kenmotsu manifold and assume that the distribution $\mathcal{D}$ has K\"{a}hler leaves. Then, the following formula holds
$$
\sum_{i=1}^{2n+1}(\nabla_{X_i}h)X_i=\varphi Q\xi
$$
where $\{X_i\}_{1\leq i\leq 2n+1}$ is an arbitrary orthonormal frame.
\end{lemma}
\begin{proof}
The differentiation of the relation $h\varphi=-\varphi h$, together with (\ref{kaehler leaves  bis}), for any vector fields $X,Y$, yields
\begin{equation*}
(\nabla_{X}h)\varphi Y+\varphi(\nabla_{X}h)Y=\eta(Y)(h\varphi X+h^2X)-\{g(\varphi X, hY)+g(hX, hY)\}\xi.
\end{equation*}
Taking $X=Y=X_{i}$, summing on $i$ and using $\mathrm{tr}(h\varphi)=0$ and $h(\xi)=0$, we get
\begin{equation}\label{tracenh0}
\sum_{i=1}^{2n+1}\{(\nabla_{X_{i}}h)\varphi X_{i}+\varphi(\nabla_{X_{i}}h)X_{i}\}=-\mathrm{tr}(h^2)\xi.
\end{equation}
Now, (\ref{tracenh'}) can be written in terms of $h$ as
$$
\sum_{i=1}^{2n+1}(\nabla_{X_{i}}h)\varphi X_{i}+\sum_{i=1}^{2n+1}h(\nabla_{X_{i}}\varphi) X_{i}=Q\xi+2\xi,
$$
from which, using (\ref{tracenfi}), we get $\sum_{i=1}^{2n+1}(\nabla_{X_{i}}h)\varphi X_{i}=Q\xi+2\xi$. Substituting this expression in (\ref{tracenh0}), we obtain
$$
\sum_{i=1}^{2n+1}\varphi(\nabla_{X_{i}}h)X_{i}=-(2+\mathrm{tr}(h^2))\xi-Q\xi.
$$
Finally, we get the required result acting by $\varphi$ and using $\sum_{i=1}^{2n+1}g((\nabla_{X_{i}}h)X_{i},\xi)=0$, which, by direct computation, follows from the fact that $g(h^2 X, \varphi X)=0$ for all vector fields $X$ and $\mathrm{tr}(h)=0$.
\end{proof}

In the three-dimensional case we have:

\begin{lemma}
If $(M^{3},\varphi,\xi,\eta,g)$ is an almost Kenmotsu manifold satisfying a generalized nullity condition, then one has
\begin{gather}
Q=aI+b\eta\otimes\xi+\mu T \label{riccitensor} \\
T(\mathrm{grad}\ \mu)=\mathrm{grad}\ \kappa-(\xi \kappa)\xi, \label{gradients}
\end{gather}
where $a=\tau/2-\kappa$ and $b=3\kappa-\tau/2$, $\tau$ being the scalar curvature of $g$, and $T$ is either $h$ or $h'$, according to which nullity condition is satisfied.
\end{lemma}
\begin{proof}
Since we know that every Riemannian $3$-manifold has vanishing Weyl conformal tensor field, we have
\begin{equation}\label{0Weyl}
\begin{split}
R(X,Y)Z&=g(Y,Z)QX-g(X,Z)QY+g(QY,Z)X\\
&\quad -g(QX,Z)Y-\frac{\tau}{2}\big(g(Y,Z)X-g(X,Z)Y\big).
\end{split}
\end{equation}
Taking into account that $Q\xi=2\kappa\xi$ (cf. (\ref{eq:h^2}) for n=1), we have
\begin{align*}
R(X,\xi)\xi=QX-2\kappa\eta(X)\xi+2kX
-2\kappa\eta(X)\xi-\frac{\tau}{2}\big(X-\eta(X)\xi\big).
\end{align*}
Comparing this expression with $R(X,\xi)\xi=\kappa(X-\eta(X)\xi)+\mu TX$, obtained by means of (\ref{nullitycond1}) if $T=h$ and (\ref{nullitycond2}) if $T=h'$, we get (\ref{riccitensor}).
Let $\{X_1=\xi, X_2=X, X_3=\varphi X\}$ be an orthonormal local frame adapted to the structure. Using (\ref{riccitensor}) and (\ref{e:nablaxi}), since $\eta(\nabla_X X+\nabla_{\varphi X} \varphi X)=-2$ and $a+b=2\kappa$, one has
\begin{align*}
\sum_{i=1}^{3}(\nabla_{X_i}Q)X_i&=\sum_{i=1}^{3}\left(X_i(a)X_i\right)+\sum_{i=1}^{3}\left(X_i(\mu)TX_i\right)+\mu\sum_{i=1}^{3}\left(\nabla_{X_i}T\right)X_i \\
& \qquad +\eta(\mathrm{grad}\ b)\xi -b\left(\eta(\nabla_X X)\right)\xi-b\left(\eta(\nabla_{\varphi X} \varphi X)\right)\xi \\
&=\mathrm{grad}\ a+T(\mathrm{grad}\ \mu)+\mu\sum_{i=1}^{3}\left(\nabla_{X_i}T\right)X_i+2\xi(\kappa)\xi-\xi(a)\xi+2b\xi.
\end{align*}
In the case of the $(\kappa,\mu)$-nullity condition $T=h$, and Lemma \ref{tracenh} implies $\sum_{i=1}^{3}\left(\nabla_{X_i}T\right)X_i=0$.
Being $a=\tau/2-\kappa$ and using the well-known formula $
(1/2)\mathrm{grad}\ \tau=\sum_{i=1}^{3}(\nabla_{X_i}Q)X_i$,
from the above equation one gets
$$
\xi(\kappa)\xi-\mathrm{grad}\ \kappa+h(\mathrm{grad}\ \mu)+\xi(\kappa)\xi-\xi(a)\xi+2b\xi=0.
$$
Since the vector field $\xi(\kappa)\xi-\mathrm{grad}\ \kappa+h(\mathrm{grad}\ \mu)$ is orthogonal to $\xi$, (\ref{gradients}) with $T=h$ follows. If $\xi$ belongs to the $(\kappa,\mu)'$-nullity distribution, then $T=h'$ and $\sum_{i=1}^{3}\left(\nabla_{X_i}T\right)X_i = 2(\kappa+1)\xi$, as it follows from (\ref{tracenh'}). Similar arguments as in previous case show that (\ref{gradients}) is still true.
\end{proof}

As a consequence of (\ref{riccitensor}), we have the following result
concerning the $\eta$-Einstein condition.

\begin{proposition}\label{equiv}
Let $(M,\varphi,\xi,\eta,g)$ be a $3$-dimensional almost Kenmotsu manifold. Then, the following conditions are equivalent:

\begin{itemize}
  \item[\rm(a)] $M$ is $\eta$-Einstein;
  \item[\rm(b)] $\xi\in N(\kappa)$ for some $\kappa\in\mathfrak{F}(M)$.
\end{itemize}
\end{proposition}

\begin{proof}
Since the generalized $(\kappa,\mu)$- and $(\kappa,\mu)'$-nullity distributions coincide with the $\kappa$-nullity distribution when $\mu=0$ or $h=0$ (equivalently $h'=0$), the implication ${\rm(b)}\Rightarrow{\rm(a)}$ immediately follows from (\ref{riccitensor}). \\
To prove ${\rm(a)}\Rightarrow{\rm(b)}$, it suffices to use (\ref{0Weyl}) and (\ref{ricci1}). From (\ref{ricci1}) we have $Q(\xi)=(a+b)\xi$ and $\tau=3a+b$. Using (\ref{0Weyl}), by direct computation
one gets:
$$
R(X,Y)\xi = \left(\frac{a+b}{2}\right)(\eta(Y)X-\eta(X)Y),
$$
for all vector fields $X,Y$, which means that $\xi$ belongs to the generalized $\kappa$-nullity distribution with $\kappa=(a+b)/2$.
\end{proof}

In \cite{TK} the concept of \emph{$N(\kappa)$-quasi Einstein manifold} is introduced as a quasi Einstein manifold with $\xi$ belonging to some $N(\kappa)$. In particular, it is proved that every $n$-dimensional conformally flat quasi Einstein manifold is a $N\left((a+b)/(n-1)\right)$-quasi Einstein manifold. Therefore, the implication ${\rm(a)}\Rightarrow{\rm(b)}$ of the above proposition can be considered as a consequence of the quoted more general result. Furthermore, we immediately deduce the following corollary of Proposition \ref{equiv}:

\begin{corollary}
Any almost Kenmotsu $3$-manifold with $\xi$ belonging to a $\kappa$-nullity distribution is a $N(\kappa)$-quasi Einstein manifold.
\end{corollary}

As proven in \cite{PS1}, a $3$-dimensional almost Kenmotsu manifold with $h=0$ always is $\eta$-Einstein. Therefore, from now on, we will restrict our investigations  mainly to the more meaningful case $h\neq0$. Note that, in this case, in the $N(\kappa)$-distribution $\kappa$ must be a non-constant function.

\section{A class of almost Kenmotsu generalized $(\kappa,\mu)$-manifolds}
\label{sec4}

In this section, we are interested in $3$-dimensional almost Kenmotsu generalized $(\kappa,\mu)$-manifolds with $\kappa$ such that $d\kappa\wedge\eta=0$. For such manifolds, the following two typical situations should be treated: either $\kappa=-1$ identically on $M$, case in which $M$ is Kenmotsu and locally a warped product of a K\"{a}hler manifold and an open interval, or $\kappa<-1$ everywhere on $M$. This follows from the following fact.

\begin{lemma}\label{property on k}
Let $(M^{3},\varphi,\xi,\eta,g)$ be an almost Kenmotsu generalized $(\kappa,\mu)$-manifold with $d\kappa\wedge\eta=0$. If $\kappa=-1$ at a certain point of $M$, then $\kappa=-1$ everywhere on $M$ and $h$ vanishes identically.
\end{lemma}
\begin{proof}
Let $Z$ be the closed subset of $M$ containing the points $q$ at which $\kappa=-1$, which is nonempty by hypothesis, and fix $p\in Z$. Then, $\lambda(p)=\sqrt{-1-\kappa(p)}=0$; furthermore $d\kappa\wedge\eta=0$ implies the same condition for $\lambda$. Therefore, choosing a coordinate neighbourhood $U$ around $p$ with coordinate $(x,y,t)$ such that $\xi=\partial/\partial t$ and $\eta=dt$, the function $\lambda$ restricted to $U$ depends only on $t$ and it satisfies the linear differential equation $d\lambda/dt=-2\lambda$. So $\lambda=ce^{-2t}$ for some real constant $c\geq0$. Since $\lambda(p)=0$, we get $c=0$, hence $\lambda=0$ and $\kappa=-1$ on the whole $U$. It follows that $U\subset Z$ and, consequently, $Z$ is also an open subset of $M$. Thus, $Z=M$ since $M$ is connected.
\end{proof}

 In view of the above considerations, we focus our attention on the case $\kappa<-1$.

\begin{proposition}\label{3-Levi-Civita}
Let $(M^{3},\varphi,\xi,\eta,g)$ be an almost Kenmotsu generalized $(\kappa,\mu)$-manifold with $h\neq0$. Then, for any unit eigenvector $X$ of $h$ with eigenvalue $\lambda=\sqrt{-1-\kappa}$, one has
\begin{itemize}
  \item[\rm(i)] $\nabla_ X \xi=X-\lambda \varphi X;\quad\ \nabla_ {\varphi X} \xi=\varphi X-\lambda X$.
  \item[\rm(ii)] $\nabla_{\varphi X} \varphi X=\frac{X(\lambda)}{2\lambda}X-\xi;\quad\ \nabla_ X X=\frac{\varphi X(\lambda)}{2\lambda}\varphi X-\xi$.
  \item[\rm(iii)] $\nabla_ X \varphi X=\lambda\xi-\frac{\varphi X(\lambda)}{2\lambda}X;\quad\ \nabla_ {\varphi X} X=\lambda\xi-\frac{X(\lambda)}{2\lambda}\varphi X$.
  \item[\rm(iv)] $\nabla_ \xi  X=-\frac{\mu}{2}\varphi X;\quad\ \nabla_ \xi \varphi X=\frac{\mu}{2} X$.
\end{itemize}
\end{proposition}

\begin{proof}
 Both eigendistributions $[\lambda]$ and $[-\lambda]$ of $h$ are $1$-dimensional and integrable; moreover, fixed a unit vector field $X \in [\lambda]$, we notice that $\{X,\varphi X, \xi\}$ is an orthonormal (local) frame for $TM$.

(i) The formulas are immediate consequences of (\ref{e:nablaxi}).

(ii) By (\ref{eq:h^2}), we have $\varphi Q\xi=0$, so that Lemma \ref{tracenh} gives
$(\nabla_X h)X+(\nabla_{\varphi X} h)\varphi X=0,$
that is
\begin{equation}\label{(*)}
X(\lambda)X+\lambda \nabla_X X-h(\nabla_X X)-\varphi X(\lambda)\varphi X-\lambda \nabla_{\varphi X} \varphi X-h(\nabla_{\varphi X} \varphi X)=0.
\end{equation}
By scalar multiplication with $X$, we get $X(\lambda)-g(\nabla_{\varphi X} \varphi X,\lambda X)-\lambda g(\nabla_{\varphi X} \varphi X, X)=0$,
from which, being $\lambda\neq0$, it follows that
$g(\nabla_{\varphi X} \varphi X, X)=X(\lambda)/(2\lambda)$.
Moreover, we have
$g(\nabla_{\varphi X} \varphi X, \xi)=-g(\varphi X, \nabla_{\varphi X}\xi)=-1.$
In this way, we obtain the first formula. As for $\nabla_X X$,
if we take the inner product of (\ref{(*)}) with $\varphi X$, we have
$\lambda g(\nabla_X X,\varphi X)-g(\nabla_X X,-\lambda\varphi X)-\varphi X(\lambda)=0$
and hence
$g(\nabla_X X,\varphi X)=\varphi X(\lambda)/(2\lambda)$.
Moreover,
$g(\nabla_X X, \xi)=-g(X,\nabla_X \xi)=-1$,
and the required formula is obtained.

(iii) The equations follow from the previous ones, considering the following relations:
\begin{align*}
&g(\nabla_X \varphi X, X)=-g(\varphi X,\nabla_X X), \quad g(\nabla_{\varphi X} X, \varphi X)=-g(X,\nabla_{\varphi X} \varphi X), \\
&g(\nabla_X \varphi X, \xi)=-g(\varphi X,\nabla_X \xi)=\lambda, \quad
g(\nabla_{\varphi X} X, \xi)=-g(X,\nabla_{\varphi X} \xi)=\lambda.
\end{align*}

(iv) Obviously we have $g(\nabla_ \xi  X,X)=0$ and $g(\nabla_ \xi  X,\xi)=0$ since $\nabla_ {\xi} \xi=0$. It remains only to compute the component along $\varphi X$. From (\ref{eq:nh}) we have
$$
\xi(\lambda)X+\lambda \nabla_\xi X-h(\nabla_\xi X)=-2\lambda X-\mu\varphi(\lambda X).
$$
The inner product with $\varphi X$ gives
$g(\lambda \nabla_\xi X,\varphi X)+g( \nabla_\xi X,\lambda\varphi X)=-\mu\lambda$,
from which, being $\lambda\neq0$, it follows that
$g( \nabla_\xi X,\varphi X)=-\mu/2$.
Finally, taking into account that $\nabla_\xi \varphi=0$, we get
$\nabla_\xi \varphi X=(\mu/2) X$.
\end{proof}

Locally, an almost Kenmotsu generalized $(\kappa,\mu)$-manifold with $d\kappa\wedge\eta=0$ and $\kappa<-1$ can be described as follows.

\begin{theorem}\label{(k,mu)-local description}
Let $(M^{3},\varphi,\xi,\eta,g)$ be an almost Kenmotsu generalized $(\kappa,\mu)$-manifold with $d\kappa\wedge\eta=0$ and $\kappa<-1$. Then, in a neighbourhood $U$ of every point $p\in M^3$ there exist coordinates $x,y,z$ with $z<-1$ and an orthonormal frame $\{\xi, X, \varphi X\}$ of eigenvectors of $h$ with $hX=\lambda X$, such that on $U$ $\kappa=z$, $\mu$ only depends on $z$ and
$$
X=\frac{\partial}{\partial x},\quad \varphi X=\frac{\partial}{\partial y}, \quad \xi=a\frac{\partial}{\partial x}+b \frac{\partial}{\partial y}-4(z+1)\frac{\partial}{\partial z},
$$
where $a=x-(1/2)(\mu+2\sqrt{-1-z})+f$ and $b=y+(1/2)(\mu-2\sqrt{-1-z})+r$, $f,r$ being smooth functions of $z$ on $U$.
\end{theorem}
\begin{proof}
Notice that $d\kappa\wedge\eta=0$ and (\ref{gradients}) imply $d\mu\wedge\eta=0$, since $h\neq0$ and $\ker(h)=\mathrm{Span}\{\xi\}$. Moreover, we have $Z(\lambda)=0$, for all $Z\in \mathcal{D}$. By Proposition \ref{3-Levi-Civita}, given an orthonormal local frame $\{\xi, X, \varphi X\}$ with $hX=\lambda X$, we get $[X,\varphi X]=-(\varphi X(\lambda)/2\lambda)X+(X(\lambda)/2\lambda)\varphi X$ which implies $[X,\varphi X]=0$. It follows that, fixed a point $p\in M^3$, there exist coordinates $(x',y',t)$ on an open neighbourhood $V$ of $p$ such that
$$
X=\frac{\partial}{\partial x'},\quad \varphi X=\frac{\partial}{\partial y'} \quad \text{and} \quad \xi=a\frac{\partial}{\partial x'}+b \frac{\partial}{\partial y'}+c\frac{\partial}{\partial t},
$$
where $a,b,c$ are smooth functions on $V$ with $c\neq0$ everywhere. Now, from the conditions $[X,\xi]\in \mathcal{D}$ and $[\varphi X,\xi]\in \mathcal{D}$ we deduce that $\partial c/\partial x'=0$ and $\partial c/\partial y'=0$. Therefore, if we consider on $V$ the linearly independent vector fields $X$, $\varphi X$ and $W:=c\partial/\partial t$, we have
$$
[X,\varphi X]=0, \quad [X,W]=0, \quad [\varphi X,W]=0.
$$
This means that there exists a coordinate system $\{U,(x,y,z')\}$ around $p$ in $V$ such that $X=\partial/\partial x$, $\varphi X=\partial/\partial y$ and $W=\partial/\partial z'$. Thus, on the open set $U$ we have
$\xi=a\partial/\partial x+b \partial/\partial y+\partial/\partial z'$. From (\ref{eq:nh}), since $\partial \lambda/\partial x=X(\lambda)=0$ and $\partial \lambda/\partial y=\varphi X(\lambda)=0$, it follows that $\lambda=c'e^{-2z'}$ and $\kappa=-1-\lambda^2=-1-\bar{c}e^{-4z'}$, for some real constants $c',\bar{c}>0$. The change of coordinates  $z=-1-\bar{c}e^{-4z'}$ gives a chart $\{U,(x,y,z)\}$ at $p$ such that $\kappa=z<-1$, $\mu=\mu(z)$ and
$$
X=\frac{\partial}{\partial x},\quad \varphi X=\frac{\partial}{\partial y} \quad \text{and} \quad \xi=a\frac{\partial}{\partial x}+b \frac{\partial}{\partial y}-4(z+1)\frac{\partial}{\partial z}.
$$
To conclude the proof, we have to calculate the functions $a,b$. To this end, we have
$$
[\xi,X]=-\frac{\partial a}{\partial x}\frac{\partial}{\partial x}-\frac{\partial b}{\partial x}\frac{\partial}{\partial y}, \qquad [\xi,\varphi X]=-\frac{\partial a}{\partial y}\frac{\partial}{\partial x}-\frac{\partial b}{\partial y}\frac{\partial}{\partial y}.
$$
On the other hand, by using Proposition \ref{3-Levi-Civita}, one has
$$
[\xi,X]=-X+\left(\lambda-\frac{\mu}{2}\right)\varphi X, \qquad [\xi,\varphi X]=\left(\lambda+\frac{\mu}{2}\right)X-\varphi X.
$$
The comparison of these relations with the previous leads to
$$
\frac{\partial a}{\partial x}=1, \qquad \frac{\partial a}{\partial y}=-\left(\lambda+\frac{\mu}{2}\right), \qquad \frac{\partial b}{\partial x}=\frac{\mu}{2}-\lambda, \qquad \frac{\partial b}{\partial y}=1.
$$
By integration of this system of differential equations, being $\lambda,\mu$ functions only depending on $z$, we get $a=x-(1/2)(\mu+2\sqrt{-1-z})+f$ and $b=y+(1/2)(\mu-2\sqrt{-1-z})+r$, for some functions $f=f(z),r=r(z)$.
\end{proof}

\begin{remark}
The above result allows us to obtain a complete local classification of $3$-dimensional almost Kenmotsu generalized $(\kappa,\mu)$-manifolds with $d\kappa\wedge\eta=0$ and $\kappa<-1$. In fact, we are going to construct in $\mathbb{R}^3$ a model for each of them as follows.

Let $M$ be the open submanifold of $\mathbb{R}^3$ defined by $M:=\{(x,y,z)\in \mathbb{R}^3|z<-1\}$ and $\mu,f,r:M \to \mathbb{R}$ be three smooth functions of $z$. Let us denote again by $x,y,z$ the coordinates induced on $M$ by the standard ones on $\mathbb{R}^3$. We consider on $M$
$$
\xi=\alpha\frac{\partial}{\partial x}+\beta \frac{\partial}{\partial y}-4(z+1)\frac{\partial}{\partial z}, \qquad \eta=-\frac{1}{4(1+z)}dz,
$$
the Riemannian metric $g$ given by
$$
g = dx\otimes dx + dy\otimes dy + (1+\alpha^2+\beta^2)\eta\otimes \eta - \alpha(dx\otimes\eta+\eta\otimes dx) - \beta(dy\otimes\eta+\eta\otimes dy)
$$
and the $(1,1)$-tensor field $\varphi$ represented, with respect to the global coordinate vector fields, by the following matrix
$$
\varphi=
\left(\begin{array}{ccc}
 0 &  -1 &  -\beta/\gamma \\
 1 &  0  &   \alpha/\gamma \\
 0 &  0  &  0
\end{array}\right),
$$
where $\alpha=x-(1/2)(\mu+2\sqrt{-1-z})+f(z)$, $\beta=y+(1/2)(\mu-2\sqrt{-1-z})+r(z)$ and $\gamma=4(1+z)$.
It is easy to check that $(M,\varphi,\xi,\eta,g)$ is an almost Kenmotsu manifold and that $\{E_1=\partial/\partial x, E_2=\partial/\partial y, E_3=\xi\}$ make up a global $\varphi$-basis on $M$, that is a $g$-orthonormal global frame such that $\varphi E_1=E_2$. Moreover, by direct computation, putting $\lambda=\sqrt{-1-z}$, we get $[E_1,E_2]=0$, $[E_1,E_3]=E_1-\left(\lambda-\mu/2\right)E_2$, $[E_2,E_3]=-\left(\lambda+\mu/2\right)E_1+E_2$, and hence $h E_1=\lambda E_1$, $h E_2=-\lambda E_2$ and $hE_3=0$.
Now, we remark that, putting $X=E_1$ and using the Koszul's formula for the Levi-Civita connection $\nabla$ of $g$, we find an orthonormal frame $\{\xi,X,\varphi X\}$ and the expressions of $\nabla$ as described in Proposition \ref{3-Levi-Civita}. Using these formulas and the definition of the curvature tensor, we finally obtain that $(M,\varphi,\xi,\eta,g)$ is a $3$-dimensional almost Kenmotsu generalized $(\kappa,\mu)$-manifold with $\kappa=z$. By virtue of Theorem \ref{(k,mu)-local description},
any $3$-dimensional almost Kenmotsu generalized $(\kappa,\mu)$-manifold is locally isomorphic to one of above manifolds.
\end{remark}

\section{A class of almost Kenmotsu generalized $(\kappa,\mu)'$-manifolds}
\label{sec5}

Here we present a local description of the structure of almost Kenmotsu generalized $(\kappa,\mu)'$-manifolds of dimension $3$ with $\kappa$ such that $d\kappa\wedge\eta=0$.

First of all, taking account of the second equation in (\ref{eq:nh'}) and arguing as in the proof of Lemma \ref{property on k}, we obtain:

\begin{lemma}
Let $(M^{3},\varphi,\xi,\eta,g)$ be an almost Kenmotsu generalized $(\kappa,\mu)'$-manifold with $d\kappa\wedge\eta=0$. If $\kappa=-1$ at a certain point of $M$, then $\kappa=-1$ everywhere on $M$ and $h'=h=0$. In particular, it follows that $M$ is locally a warped product of an open interval and a K\"{a}hler manifold.
\end{lemma}

Therefore, also for almost Kenmotsu generalized $(\kappa,\mu)'$-manifolds we shall discuss the case $\kappa<-1$, or equivalently $h'\neq0$ everywhere. We first obtain the following proposition.

\begin{proposition}\label{3-Levi-Civita'}
Let $(M^{3},\varphi,\xi,\eta,g)$ be an almost Kenmotsu generalized $(\kappa,\mu)'$-manifold with $h'\neq0$. Then, for any unit $X\in [\lambda]'$, $\lambda=\sqrt{-1-\kappa}$, one has
\begin{itemize}
\item[(i)] $\nabla_ X \xi=(1+\lambda)X;\quad\ \nabla_ {\varphi X} \xi=(1-\lambda)\varphi X$
\item[\rm(ii)] $\nabla_ X \varphi X=-\frac{\varphi X(\lambda)}{2\lambda}X;\quad\ \nabla_ {\varphi X} X=-\frac{ X(\lambda)}{2\lambda}\varphi X$.
\item[(iii)] $\nabla_{\varphi X} \varphi X=\frac{ X(\lambda)}{2\lambda}X-(1-\lambda)\xi;\quad\ \nabla_ X X=\frac{\varphi  X(\lambda)}{2\lambda}\varphi X-(1+\lambda)\xi$.
\item[(iv)] $\nabla_ \xi  X=0;\quad\ \nabla_ \xi \varphi X=0$.
\end{itemize}
\end{proposition}

\begin{proof}
Let $X$ be a unit eigenvector of $h'$ corresponding to the eigenvalue $\lambda$. Then, since the distributions $[\lambda]'$ and $[-\lambda]'$ are $1$-dimensional, $\{X,\varphi X, \xi\}$ is an orthonormal (local) frame.

(i) This is a direct consequence of (\ref{e:nablaxi}).

(ii) Obviously $g(\nabla_ X \varphi X,\varphi X)=0$ and $g(\nabla_ X \varphi X,\xi)=-g(\varphi X,(1+\lambda)X)=0$, so that $\nabla_ X \varphi X\in \mathrm{Span}\{X\}$. Equation (\ref{eq:h' codazzi}) with $Y=\varphi X$ implies
\begin{equation}\label{(**)}
-X(\lambda)\varphi X-\lambda \nabla_X \varphi X-h'(\nabla_X \varphi X)-\varphi X(\lambda) X-\lambda \nabla_{\varphi X} X+h'(\nabla_{\varphi X} X)=0.
\end{equation}
By inner product with $X$, we get $-g(\lambda \nabla_X \varphi X, X)-g( \nabla_X \varphi X,\lambda X)-\varphi X(\lambda)=0$. Therefore, being $\lambda\neq0$, we find $g(\nabla_X \varphi X, X)=-\varphi X(\lambda)/2\lambda$. Analogously, $g(\nabla_ {\varphi X}  X,X)=0$, $g(\nabla_ {\varphi X}  X, \xi)=0$ and the scalar product of (\ref{(**)}) with $\varphi X$ imply the second formula.

(iii) We get this using the previous equations. Indeed, we have
$$
\begin{array}{l}
g(\nabla_{\varphi X} \varphi X,X)=-g(\varphi X,\nabla_{\varphi X} X)=\frac{ X(\lambda)}{2\lambda}, \quad g(\nabla_X X,\varphi X)=-g( X,\nabla_X \varphi X)=\frac{\varphi X(\lambda)}{2\lambda}, \\
\\
g(\nabla_{\varphi X} \varphi X, \xi)=-g(\varphi X,\nabla_{\varphi X} \xi)=-(1-\lambda), \quad
g(\nabla_X X, \xi)=-g(X,\nabla_X \xi)=-(1+\lambda).
\end{array}
$$

(iv) Since $g(\nabla_\xi X, \xi)=0$, we have $\nabla_\xi X \in \mathrm{Span}\{\varphi X\}$. Using the first equation in (\ref{eq:nh'}), we obtain
$$
\xi(\lambda)X+\lambda\nabla_\xi X-h'(\nabla_\xi X)=-\lambda(\mu+2)X.
$$
Taking the scalar product with $\varphi X$, we get $2\lambda g(\nabla_\xi X,\varphi X)=0$. Being $\lambda\neq0$, we have $\nabla_\xi X=0$. From $\nabla_\xi \varphi=0$, the remaining relation follows.
\end{proof}

\begin{theorem}\label{(k,mu)'-local description}
Let $(M^{3},\varphi,\xi,\eta,g)$ be an almost Kenmotsu generalized $(\kappa,\mu)'$-manifold. If $d\kappa\wedge\eta=0$ and $\kappa<-1$, then one has:
\begin{itemize}
\item[(i)] If $\mu=-2$, then $\kappa$ is constant.
\item[(ii)] If $\mu\neq-2$, then in a neighbourhood $U$ of every point $p\in M^3$ there exist coordinates $x,y,z$ with $z<-1$ and an orthonormal frame $\{\xi, X, \varphi X\}$ of eigenvectors of $h'$ with $h'X=\lambda X$, such that $\kappa=z$, $\mu$ only depends on $z$ on $U$ and
$$
X=\frac{\partial}{\partial x},\quad \varphi X=\frac{\partial}{\partial y}, \quad \xi=a\frac{\partial}{\partial x}+b \frac{\partial}{\partial y}-2(z+1)(\mu+2)\frac{\partial}{\partial z},
$$
where $a=x(1+\sqrt{-1-z})+f$ and $b=y(1-\sqrt{-1-z})+r$, $f,r$ being smooth functions of $z$ on $U$.
\end{itemize}
\end{theorem}
\begin{proof}
Since $\lambda^2=-1-\kappa$, the hypothesis $d\kappa\wedge\eta=0$ implies $d\lambda\wedge\eta=0$, or equivalently $Z(\lambda)=0$ for all $Z\in \mathcal{D}$. From (\ref{eq:nh'}) we deduce $d\mu\wedge\eta=0$. If $\mu=-2$, by (\ref{eq:nh'}), we have $\xi(\lambda)=0$ which implies $\lambda$ constant and so $\kappa$ constant.

Now, we assume $\mu\neq-2$. By Proposition \ref{3-Levi-Civita'}, given an orthonormal local frame $\{\xi, X, \varphi X\}$ with $h'X=\lambda X$, we get $[X,\varphi X]=-(\varphi X(\lambda)/2\lambda)X+(X(\lambda)/2\lambda)\varphi X$ which implies $[X,\varphi X]=0$. Consequently, fixed a point $p\in M^3$, there exist coordinates $(\bar{x},\bar{y},\bar{t})$ on an open neighbourhood $V$ of $p$ such that
$$
X=\frac{\partial}{\partial \bar{x}},\quad \varphi X=\frac{\partial}{\partial \bar{y}} \quad \text{and} \quad \xi=a\frac{\partial}{\partial \bar{x}}+b \frac{\partial}{\partial \bar{y}}+c\frac{\partial}{\partial \bar{t}},
$$
where $a,b,c$ are smooth functions on $V$ with $c\neq0$ everywhere on $V$. From the conditions $[X,\xi]\in \mathcal{D}$ and $[\varphi X,\xi]\in \mathcal{D}$ we deduce that $\partial c/\partial \bar{x}=0$ and $\partial c/\partial \bar{y}=0$. Therefore, the Lie brackets of the vector field $W:=c\partial/\partial \bar{t}$ with the other two coordinate vector fields vanish, thus obtaining coordinates $(x',y',t')$ on an open neighbourhood $U'$ of $p$ in $V$ such that $X=\partial/\partial x'$, $\varphi X=\partial/\partial y'$ and $W=\partial/\partial t'$. Thus on $U'$ we have
$\xi=a\partial/\partial x'+b \partial/\partial y'+\partial/\partial t'$.
Moreover, $\lambda|_{U'}$ depends only on $t'$ and the second equation in (\ref{eq:nh'}) reads
\begin{equation}\label{diffeq2}
\frac{d \lambda}{d t'}=-\lambda(\mu+2),
\end{equation}
since $\partial \lambda/\partial x'=0$ and $\partial \lambda/\partial y'=0$. Next, considering the vector field $W':=1/(\mu+2)\partial/\partial t'$, since $\mu$ independent of $x'$ and $y'$, we have $[X,W']=0$ and $[\varphi X,W']=0$. It follows that it is possible to find a chart $\{U,(x,y,z')\}$ at $p$ in $U'$ such that $X=\partial/\partial x$, $\varphi X=\partial/\partial y$, $W'=\partial/\partial z'$ and $\xi=a\partial/\partial x + b \partial/\partial y+(\mu+2)\partial/\partial z'$, again denoting by $a,b$ the restriction to $U$ of these functions. With respect to these coordinates, (\ref{diffeq2}) becomes $d \lambda/d z'=-\lambda$, from which we get $\lambda=c'e^{-z'}$, with $c'>0$ a real constant. Hence $\kappa=-1-\bar{c}e^{-2z'}$ for some real constant $\bar{c}>0$. Finally, the substitution $z=-1-\bar{c}e^{-2z'}$ gives the desired chart $\{U,(x,y,z)\}$ at $p$. To conclude the proof, it remains to compute the functions $a,b$. To do this, let us explicit the Lie brackets $[\xi,\varphi X]$ and $[\xi,X]$, obtaining
$$
[\xi,X]=-\frac{\partial a}{\partial x}\frac{\partial}{\partial x}-\frac{\partial b}{\partial x}\frac{\partial}{\partial y}, \qquad [\xi,\varphi X]=-\frac{\partial a}{\partial y}\frac{\partial}{\partial x}-\frac{\partial b}{\partial y}\frac{\partial}{\partial y}.
$$
On the other hand, by using Proposition \ref{3-Levi-Civita'}, we have
$$
[\xi,X]=-(1+\lambda)X, \qquad [\xi,\varphi X]=-(1-\lambda)\varphi X.
$$
Comparing these relations with the previous, we get
$$
\frac{\partial a}{\partial x}=1+\lambda, \qquad \frac{\partial a}{\partial y}=0, \qquad \frac{\partial b}{\partial x}=0, \qquad \frac{\partial b}{\partial y}=1-\lambda.
$$
The integration of this system yields $a=x(1+\sqrt{-1-z})+f$ and $b=y(1-\sqrt{-1-z})+r$, where $f,r$ are arbitrary smooth functions of $z$ on $U$.
\end{proof}

\begin{remark}
We may build local models for each $3$-dimensional almost Kenmotsu generalized $(\kappa,\mu)'$-manifold with $d\kappa\wedge\eta=0$, $\kappa<-1$ and $\mu\neq-2$. Let us consider  $M^3=\{(x,y,z)\in \mathbb{R}^3|z<-1\}$ and smooth functions $\mu,f,r:M \to \mathbb{R}$ depending on $z$ such that $\mu\neq-2$. We take the following vector fields
$$
E_1:=\frac{\partial}{\partial x},\quad E_2:=\frac{\partial}{\partial y}, \quad E_3:=a\frac{\partial}{\partial x}+b \frac{\partial}{\partial y}-2(z+1)(\mu+2)\frac{\partial}{\partial z},
$$
where $a=x(1+\sqrt{-1-z})+f(z)$ and $b=y(1-\sqrt{-1-z})+r(z)$.
Setting $\lambda=\sqrt{-1-z}$, we have
\begin{equation}\label{brackets}
[E_1,E_2]=0, \quad [E_1,E_3]=(1+\lambda)E_1, \quad [E_2,E_3]=(1-\lambda)E_2.
\end{equation}
Let $g$ be the Riemannian metric on $M^3$ that makes the basis $\{E_1,E_2,E_3\}$ orthonormal.
The structure tensor fields $\varphi,\xi,\eta$ are defined by putting
$$
\xi=E_3, \qquad \eta=-\frac{1}{2(1+z)(\mu+2)}dz
$$
$$
\varphi(\xi)=0, \quad \varphi(E_1)=E_2, \quad \varphi(E_2)=-E_1
$$
Then, as one can easily prove, $(M,\varphi,\xi,\eta,g)$ is in fact an almost Kenmotsu manifold. Furthermore, computing the tensor field $h'$, we get $h'E_1=\lambda E_1$ and $h'E_2=-\lambda E_2$.
From this it follows that the Levi-Civita connection of $g$, computed by means of the Koszul's formula and (\ref{brackets}), satisfies the equations stated in Proposition \ref{3-Levi-Civita'} with $X=E_1$. Finally, the direct computation of the curvature tensor shows that $(M^3,\varphi,\xi,\eta,g)$ is an almost Kenmotsu generalized $(\kappa,\mu)'$-manifold with $\kappa=z$ and $d\kappa\wedge\eta=0$.
\end{remark}

\section{Further characterizations and local models}\label{sec6}

This section is devoted to obtaining, using Darboux-like coordinates, another explicit local description of $3$-dimensional almost Kenmotsu generalized $(\kappa,\mu)$- and $(\kappa,\mu)'$-manifolds with $d\kappa\wedge\eta=0$. To this purpose, we need the following useful property.

\begin{proposition}\label{flat leaves}
Let $(M^{3},\varphi,\xi,\eta,g)$ be an almost Kenmotsu generalized $(\kappa,\mu)$- or $(\kappa,\mu)'$-manifold with $d\kappa\wedge\eta=0$ and $\kappa<-1$. Then the leaves of the canonical foliation of $M$ are flat K\"{a}hlerian manifolds.
\end{proposition}
\begin{proof}
Let $M'$ be a leaf of $\mathcal{D}$ and $(J,G)$ be the induced almost Hermitian structure.
Being almost K\"{a}hler of dimension $2$, $M'$ is obviously a K\"{a}hler manifold. In order to prove the flatness of $(M',G)$, we consider the Weingarten operator $A$ of $M'$ given by $AX=-X+\varphi h X=-X-h'X$. For a unit vector field $X\in[\lambda]$
(or $X\in[\lambda]'$), by using the Gauss equation, the sectional curvature $K'$ of $G$ is given by
$$
K'(X,\varphi X) = K(X,\varphi X) + 1 - \lambda^2.
$$
Now, using Proposition \ref{3-Levi-Civita} for the case of an almost Kenmotsu generalized $(\kappa,\mu)$-manifold and Proposition \ref{3-Levi-Civita'} for the case of an almost Kenmotsu generalized $(\kappa,\mu)'$-manifold, since $X(\lambda)=\varphi X(\lambda)=0$, one gets
$R(X,\varphi X)\varphi X = - (1 - \lambda^2)X$.
So we obtain $$K(X,\varphi X)=g(R(X,\varphi X)\varphi X, X)=-(1-\lambda^2)$$ which implies $K'(X,\varphi X)=0$.
\end{proof}

Let us now state the following characterization.

\begin{theorem}\label{local characterization}
Let $(M^{3},\varphi,\xi,\eta,g)$ be an almost contact metric manifold with $h\neq0$, and $\kappa,\mu\in\mathfrak{F}(M)$ such that $d\kappa\wedge\eta=0$.
Then, $M^{3}$ is an almost Kenmotsu generalized $(\kappa,\mu)$-manifold if and only if for any point $p\in M$, there exists an open neighbourhood $U$ of $p$ with coordinates $x_1,x_2,t$ such that $\kappa=-1-e^{-4t}$, $\mu$ only depends on $t$  and the tensor fields of the structure are expressed in the following way
\begin{equation}\label{local structure1}
\varphi=\sum_{i,j=1}^{2}\varphi^{i}_{j}dx_{j}\otimes\frac{\partial}{\partial x_{i}}, \quad \xi=\frac{\partial}{\partial t}, \quad \eta=dt, \quad g=dt\otimes dt+\sum_{i,j=1}^{2} g_{ij}dx_{i}\otimes dx_{j},
\end{equation}
where $\varphi^{i}_{j},g_{ij}$ are functions only of $t$; the fundamental $2$-form $\Phi$ is given by
\begin{equation}\label{fundform}
\Phi=e^{2t} dx_{1} \wedge dx_{2},
\end{equation}
and the non-zero components $h^{i}_{j}, B^{i}_{j}$ in $U$ of $h$ and $B:=\varphi h$, respectively, are functions of $t$ satisfying the condition $\sum_s B^{i}_{s}B^{s}_{j}=e^{-4t}\delta_{j}^{i}$ and the following system of differential equations
\begin{equation}\label{fi,h,B}
\frac{d \varphi^{i}_{j}}{dt}=2h^{i}_{j}, \quad \frac{d h^{i}_{j}}{dt}=2\lambda^2\varphi^{i}_{j}-2h^{i}_{j}-\mu B^{i}_{j}, \quad \frac{d B^{i}_{j}}{dt}=\mu h^{i}_{j}-2B^{i}_{j},
\end{equation}
where $\lambda=e^{-2t}$.
\end{theorem}
\begin{proof}
Let $M$ be an almost Kenmotsu generalized $(\kappa,\mu)$-manifold. Owing to the integrabiliy of $\mathcal{D}$ and $[\xi]$, the decomposition $TM=\mathcal{D}\oplus[\xi]$ implies that any point of $M$ admits a coordinate neighbourhood $U$ of the form $U'\times \ ]-\varepsilon,\varepsilon[$ with coordinates $x_1,x_2,t$, where $x_1,x_2$ are cooordinate on $U'$ and $t$ on $]-\varepsilon,\varepsilon[$, such that $\xi=\partial/\partial t$ and $\eta=dt$. With respect to these coordinates, the shape of $\varphi$ and $g$ is as in (\ref{local structure1}), where $\varphi^{i}_{j}$ and $g_{ij}$ are functions of $x_1,x_2,t$ in general. Denoting by $h^{i}_{j}, B^{i}_{j}$ the non-zero components in $U$ of $h$ and $B$, from (\ref{Lie derivatives1}) we get
\begin{equation}\label{differential equations}
\frac{\partial \varphi^{i}_{j}}{\partial t}=2h^{i}_{j}, \quad \frac{\partial h^{i}_{j}}{\partial t}=2\lambda^2\varphi^{i}_{j}-2h^{i}_{j}-\mu B^{i}_{j}, \quad \frac{\partial B^{i}_{j}}{\partial t}=\mu h^{i}_{j}-2B^{i}_{j}.
\end{equation}
Now, let $t_{0}\in \ ]-\varepsilon,\varepsilon[$  and consider the subset $U'\times\{t_{0}\}\subset U$ which is an open submanifold of a leaf of $\mathcal{D}$. Then, the induced complex structure $J$ has components $\varphi^{i}_{j}(t_{0}, \cdot)$ and, by virtue of Proposition \ref{flat leaves}, it can be assumed that $x_1,x_2$ are chosen in such a way that on $U'\times\{t_{0}\}$ one has $J(\partial/\partial x_1)=\partial/\partial x_2$, $J(\partial/\partial x_2)=-\partial/\partial x_1$ and the induced metric $G$ has constant components. This implies that the $\varphi_j^i$'s and $g_{ij}$'s depend on $t$ alone. Consequently, by (\ref{differential equations}), even $h_j^i$ and $B_j^i$ are functions only of $t$. Since $d\kappa\wedge\eta=0$ and $\lambda^2=-1-\kappa$, we have $\partial \lambda/\partial x_1=\partial \lambda/\partial x_2=0$, so from (\ref{eq:nh}) we get $\lambda=e^{-2t}$. Hence $\kappa=-1-e^{-4t}$. Moreover, from (\ref{gradients}) we deduce that also $\mu$ only depends on $t$. The components of $\Phi$ are all zero, except for $\Phi_{12}$. From (\ref{almKen}) one gets $\Phi_{12}=ce^{2t}$, for some real constant $c\neq0$. Up to change $x_1,x_2$ with $x'_1=\sqrt{|c|}x_1$ and $x'_2=\sqrt{|c|}x_2$, we can take $c=1$. In this way, we get the desired chart around $p$. Finally, (\ref{differential equations}) gives (\ref{fi,h,B}) and (\ref{eq:h^2}) gives $\sum_s B^{i}_{s}B^{s}_{j}=e^{-4t}\delta_{j}^{i}$. \\
Conversely, suppose that $M$ carries a structure locally represented as in (\ref{local structure1})-(\ref{fi,h,B}).
Obviously $d\eta=0$, while $d\Phi=2\eta\wedge \Phi$ follows from (\ref{fundform}).
Now, we show that $M$ satisfies the generalized $(\kappa,\mu)$-condition. In order to do this, we notice that $X_1=\partial/\partial x_1$ and $X_2=\partial/\partial x_2$ are Killing vector fields.
Hence $g(\nabla_{X_i}X_j,X_q)=0$, for any $i,j,q\in\{1,2\}$. Since the distribution orthogonal to $\xi=\partial/\partial t$ is spanned by $X_1$ and $X_2$, it follows that $\nabla_{X_i}X_j \in [\xi]$ for all $i,j\in\{1,2\}$, so that $B(\nabla_{X_i}X_j)=0$. Consequently, for the Levi-Civita connection $\nabla$ determined by $g$, we have
$$
\nabla_{X_i}X_j=\nabla_{X_j}X_i=-g(X_i,X_j-BX_j)\xi, \quad \nabla_{\xi}X_i=\nabla_{X_i}\xi=X_i-BX_i, \quad \nabla_{\xi}\xi=0.
$$
By (\ref{curv1}), we compute
\begin{align*}
R(X_i,X_j)\xi&=(\nabla_{X_j}B)X_i-(\nabla_{X_i}B)X_j=\nabla_{X_j}BX_i-\nabla_{X_i}BX_j \\
             &=\sum_{q=1}^{2}\Big(B^q_i\nabla_{X_j}X_q-B^q_j\nabla_{X_i}X_q\Big) \\
             &=-\sum_{q=1}^{2}\Big(B^q_ig(X_j,X_q-BX_q)-B^q_jg(X_i,X_q-BX_q)\Big)\xi \\
             &=-g(X_j,BX_i-B^2X_i)\xi+g(X_i,BX_j-B^2X_j)\xi=0
\end{align*}
and, using (\ref{fi,h,B}),
\begin{align*}
R(X_i,\xi)\xi&=-\nabla_{\xi}\nabla_{X_i}\xi=-X_i+BX_i+\sum_{q=1}^{2}\left(\frac{dB^q_i}{dt}X_q+B^q_i\nabla_{\xi}X_q\right) \\
             &=-X_i+BX_i+\mu hX_i-2BX_i+BX_i-B^2X_i \\
             &=(-1-\lambda^2)X_i+\mu hX_i=\kappa X_i+\mu hX_i,
\end{align*}
thus concluding the proof.
\end{proof}

Following a method analogous to that used by Dacko and Olszak in \cite{DO2} for almost cosymplectic structures, we shall construct a model as follows.

Let us consider the following three constant matrices
$$
M_1=
\begin{pmatrix}
 1 & 0  \\
 0 & -1
\end{pmatrix}, \quad
M_2=
\begin{pmatrix}
 0  & 1  \\
 -1 & 0
\end{pmatrix}, \quad
M_3=
\begin{pmatrix}
 0 & 1  \\
 1 & 0
\end{pmatrix},
$$
a smooth function $\bar{\mu}:\mathbb{R} \to \mathbb{R}$ defined on the real line with coordinate $t$, and the function $\bar{\lambda}=e^{-2t}$. Then, we define three functional matrices $F,H,B$ of order $2$
$$
F(t)=[\varphi_j^i(t)], \quad H(t)=[h_j^i(t)], \quad B(t)=[b_j^i(t)]
$$
in such a way they satisfy the following system of linear differential equations with the given initial conditions
\begin{equation}\label{differential system}
\begin{alignedat}{3}
&F'=2H,  \quad\quad & &H'=2\bar{\lambda}^2F-2H-\bar{\mu} B, \quad\quad & &B'=\bar{\mu} H-2B, \\
&F(0)=M_2,      & &H(0)=-M_3,                      & &B(0)=M_1.
\end{alignedat}
\end{equation}
Here $A'$ denotes the matrix whose coefficients are the derivatives with respect to $t$ of the entries of the matrix $A$.

Now, let $(f_i,h_i,b_i)_{1 \leq i \leq 3}$ be the unique solution defined on an open interval $]a,b[$, $a,b\in\mathbb{R}$, containing $0$, of the system of linear differential equations
$$
\left\{ \begin{array}{l}
f'_i=2h_i \\
h'_i=2\bar{\lambda}^2f_i-2h_i-\bar{\mu} b_i  \quad i=1,2,3 \\
b'_i=\bar{\mu} h_i-2b_i \end{array}
\right. \,,
$$
satisfying the following initial conditions
\begin{equation*}
\begin{alignedat}{3}
f_1(0)&=0, \quad \quad & h_1(0)&=0, \quad\quad & b_1(0)&=1, \\
f_2(0)&=1, & h_2(0)&=0, & b_2(0)&=0, \\
f_3(0)&=0, & h_3(0)&=-1,& b_3(0)&=0.
\end{alignedat}
\end{equation*}
It can be easily seen that
\begin{equation}\label{solution}
\begin{split}
F(t)&=f_1(t)M_1+f_2(t)M_2+f_3(t)M_3 \\
H(t)&=h_1(t)M_1+h_2(t)M_2+h_3(t)M_3 \\
B(t)&=b_1(t)M_1+b_2(t)M_2+b_3(t)M_3,
\end{split}
\end{equation}
is the unique solution on $]a,b[$ of (\ref{differential system}).

Moreover, let $G$ be the matrix defined on $]a,b[$ by
\begin{equation}\label{metric}
G=-M_{2}F.
\end{equation}

\begin{lemma}\label{metric matrix}
The matrices $F,H,B$ are linked by the following algebraic relations
\begin{equation}\label{algebraic relations}
\begin{alignedat}{3}
&F^2=-I_2, \quad\quad & &H^2=\bar{\lambda}^2I_2, \quad \quad & &B^2=\bar{\lambda}^2I_2, \\
&HF=-FH,              & &BF=-FB,                       & &BH=-HB,  \\
&F=\bar{\lambda}^{-2}BH,    & &H=BF,                         & &B=FH,
\end{alignedat}
\end{equation}
where $I_2$ denotes the unit real matrix. Moreover, the matrix $G(t)$ is symmetric and positive definite for any $t\in\ ]a,b[$. Explicitly, $G$ is given by
\begin{equation}\label{metric bis}
G=
\begin{pmatrix}
 f_2-f_3 & f_1   \\
 f_1     &   f_2+f_3
\end{pmatrix}.
\end{equation}
\end{lemma}

\begin{proof}
We consider the following auxiliary matrices:
\begin{equation*}
\begin{alignedat}{3}
&X_1=F^2,\quad\quad & &X_2=H^2\quad \quad & &X_3=B^2, \\
&X_4=HF+FH,           & &X_5=BF+FB,            & & X_6=BH+HB, \\
&X_7=BH-\bar{\lambda}^{2}F, \quad  & &X_8=BF-H,\quad           & &X_9=FH-B.
\end{alignedat}
\end{equation*}
Using (\ref{differential system}), it is easy to check that the following homogeneous system of linear differential equations is satisfied
\begin{equation}\label{aux system}
\left\{
\begin{array}{l}
X'_1=2X_4, \quad  X'_2=2\bar{\lambda}^{2}X_4-4X_2-\bar{\mu} X_6, \quad X'_3=\bar{\mu} X_6-4 X_3 \\
X'_4=4\bar{\lambda}^{2}X_1+4X_2-2X_4-\bar{\mu} X_5, \quad X'_5=\bar{\mu} X_4-2X_5+2X_6 \\
X'_6=2\bar{\mu} X_2-2\bar{\mu} X_3+2\bar{\lambda}^{2}X_5-4X_6,  \\
X'_7=\bar{\mu} X_2-\bar{\mu} X_3-4X_7+2\bar{\lambda}^{2}X_8,  \\
X'_8=2\bar{\lambda}^{2}X_7-2X_8+\bar{\mu} X_4-\bar{\mu} X_9, \\ X'_9=2\bar{\lambda}^{2}X_1+2X_2+\bar{\mu} X_5-\bar{\mu} X_8-2X_9.
\end{array} \right.
\end{equation}
with the initial conditions
\begin{gather*}
X_1(0)=-I_2, \quad X_2(0)=I_2, \quad X_3(0)=I_2 \\
X_4(0)=X_5(0)=X_6(0)=X_7(0)=X_8(0)=X_9(0)=O_2.
\end{gather*}
Since this system admits a unique solution, a simple verification shows that
\begin{gather*}
X_1(t)=-I_2, \quad X_2(t)=\bar{\lambda}^{2}I_2, \quad X_3(t)=\bar{\lambda}^{2}I_2 \\
X_4(t)=X_5(t)=X_6(t)=X_7(t)=X_8(t)=X_9(t)=O_2,
\end{gather*}
is the solution of (\ref{aux system}) for any $t\in\ ]a,b[$, which gives (\ref{algebraic relations}).

Applying (\ref{solution}), we obtain $G=f_1M_3+f_2I_2-f_3M_1,$
which is equivalent to the form presented in (\ref{metric bis}). Therefore, we immediately see that $G(t)$ is symmetric for any $t\in\ ]a,b[$.
Computing $F^2$ with (\ref{solution}), we find
$
F^2=(f_1^2+f_3^2-f_2^2)I_2.
$
On the other hand $F^2=-I_2$, so that $(f_2+f_3)(f_2-f_3)=1+f_1^2$. Since $f_2(0)=1$ and $f_3(0)=0$, it follows that $f_2-f_3>0$ at any point of $]a,b[$ and so $G(t)$ is positive definite for any $t\in\ ]a,b[$, completing the proof.
\end{proof}

Let now $M=\ ]a,b[\ \times \ \mathbb{R}^2\subset\mathbb{R}^3$, and denote by $(t,x_{1},x_{2})$ the coordinate global system induced on $M$ by the canonical one on $\mathbb{R}^3$.
We introduce on $M$ a structure $(\varphi,\xi,\eta,g)$
by setting:
$$
\xi:=\frac{\partial}{\partial t}, \quad \eta:=dt, \quad g:=dt\otimes dt+e^{2t}\sum_{i,j=1}^{2} G_{ij}dx_{i}\otimes dx_{j},
$$
where the $G_{ij}$'s are the coefficients of the  matrix $G$,
and $\varphi$ represented, with respect to the frame $\Big\{\frac{\partial}{\partial t},\frac{\partial}{\partial x_1},\frac{\partial}{\partial x_2}\Big\}$, by the matrix
$$
\varphi=
\begin{pmatrix}
 0 &  0 &  0 \\
 0 & \varphi_1^1  & \varphi_2^1 \\
 0 &  \varphi_1^2  &  \varphi_2^2
\end{pmatrix},
$$
$\varphi_j^i$'s being the coefficients of $F$.
From Lemma \ref{metric matrix} we have that $g$ defines a Riemannian metric tensor on $M$. Furthermore, consider the smooth functions $\mu,\lambda:M\to\mathbb{R}$ defined by
$\mu(t,x,y)=\bar{\mu}(t)$ and $\lambda(t,x,y)=\bar{\lambda}(t).$

\begin{proposition}
$(M, \varphi,\xi,\eta,g)$ is an almost Kenmotsu generalized $(\kappa,\mu)$-manifold with $\kappa=-1-e^{-4t}$ and $d\kappa\wedge\eta=0$.
\end{proposition}
\begin{proof}
Routine computations show that $(\varphi,\xi,\eta,g)$ is an almost contact metric structure on $M$.
We shall prove that it is an almost Kenmotsu generalized $(\kappa,\mu)$-structure by using Theorem \ref{local characterization}. Let us choose the neighbourhood $U=M$ for any point of $M$. By (\ref{metric}) and (\ref{algebraic relations}), one has
$$
g\left(\frac{\partial}{\partial x_1},\varphi\frac{\partial}{\partial x_2}\right)=e^{2t}\sum_{r=1}^{2}\varphi_2^rg_{1r}=e^{2t}(GF)_2^1=e^{2t},
$$
so that the fundamental $2$-form $\Phi$ satisfies (\ref{fundform}).
Computing the tensor field $h$, we see that its components, with respect to the fixed coordinates, coincide exactly with the coefficients of the matrix $H$. Furthermore, by virtue of $B=FH$ (cf. (\ref{algebraic relations})), the components of the tensor field $\varphi h$ are just the coefficients of the matrix $B$. Therefore, since the matrices $F,H,B$ satisfy (\ref{differential system}), the tensor fields $\varphi, h, \varphi h$ fulfill (\ref{fi,h,B}). Consequently, by Theorem \ref{local characterization}, the considered structure is an almost Kenmotsu generalized  $(-1-e^{-4t},\mu)$-structure on $M$.
\end{proof}

For the case of almost Kenmotsu generalized $(\kappa,\mu)'$-manifolds, we suppose $\mu\neq-2$, since if $\mu=-2$ then from (\ref{eq:nh'}) we obtain $\kappa$ constant and Theorem \ref{classification theorem} applies. By using (\ref{Lie derivatives2}) and Proposition \ref{flat leaves}, an analogue of Theorem \ref{local characterization} can be proved with similar arguments. Thus, omitting the proof, we have

\begin{theorem}
Let $(M^{3},\varphi,\xi,\eta,g)$ be an almost contact metric manifold with $h\neq0$, and $\kappa,\mu\in\mathfrak{F}(M)$ such that $d\kappa\wedge\eta=0$ and $\mu\neq-2$.
Then, $M$ is an almost Kenmotsu generalized $(\kappa,\mu)'$-manifold if and only if for any point $p\in M$, there exists an open neighbourhood $U$ of $p$ with coordinates $x_{1},x_{2},t$ such that $\mu$ only depends on $t$, $\kappa=-1-e^{-2f}$, where $f$
satisfies the equation $df=(\mu+2)dt$,
and on $U$ the tensor fields of the structure can be written as
\begin{equation*}
\varphi=\sum_{i,j=1}^{2}\varphi^{i}_{j}dx_{j}\otimes\frac{\partial}{\partial x_{i}}, \quad \xi=\frac{\partial}{\partial t}, \quad \eta=dt, \quad g=dt\otimes dt+\sum_{i,j=1}^{2} g_{ij}dx_{i}\otimes dx_{j},
\end{equation*}
where $\varphi^{i}_{j},g_{ij}$ are functions only of $t$; $\Phi$ is given by
$$
\Phi=e^{2t} dx_1 \wedge dx_2,
$$
and the non-zero components  $h^{i}_{j}, B^{i}_{j}$ in $U$ of $h$ and $h'$, respectively, are functions of $t$ satisfying the condition $\sum_s B^{i}_{s}B^{s}_{j}=e^{-2f}\delta_{j}^{i}$ and the following system of differential equations
\begin{equation}\label{fi,h',h}
\frac{d \varphi^{i}_{j}}{dt}=2h^{i}_{j}, \quad \frac{d B^{i}_{j}}{dt}=-(\mu+2) B^{i}_{j}, \quad \frac{d h^{i}_{j}}{dt}=2\lambda^2 \varphi^{i}_{j}-(\mu+2)h^{i}_{j},
\end{equation}
where $\lambda=e^{-f}$.
\end{theorem}

\begin{remarks}
\begin{enumerate}[1.]
\item As a particular case, we recover Example 6.3 from \cite{PS2}.
\item Clearly, a construction similar to that carried for almost Kenmotsu generalized $(\kappa,\mu)$-manifolds can also be made in this case, starting with the system of differential equations corresponding to (\ref{fi,h',h}).
\end{enumerate}
\end{remarks}

\proof[Acknowledgements]
The author would like to express his gratitude to Prof. A. M. Pastore
for her valuable comments and helpful suggestions for improving the paper.
He also wish to thank the anonymous referee for the useful improvements suggested.

\end{document}